\documentclass[twocolumn]{svjour3}

\smartqed

\usepackage{fix-cm}
\usepackage{graphicx}
\usepackage{fixltx2e}
\usepackage{dblfloatfix}
\usepackage{subfig}
\usepackage{tabularx}
\usepackage{natbib}             % Use citep, citet in citations
\usepackage{amsmath,mathtools}  % Cases
\usepackage{color}              % Color text to indicate changes
\usepackage{multirow}           % More complex tables
\usepackage{chngcntr}           % Show section in table

\counterwithout{table}{section}

\journalname{}

\begin{document}

%%%%%%%%%%%%%%%%%%%%%%%%%%%%%%%%%%%%%%%%%%%%%%%%%%%%%%%%%%%%%%%%%%%
% Title
%%%%%%%%%%%%%%%%%%%%%%%%%%%%%%%%%%%%%%%%%%%%%%%%%%%%%%%%%%%%%%%%%%%

\title{Density and Level Set-XFEM Schemes for Topology Optimization of 3-D Structures}

%\subtitle{Topology Optimization in 3D}

\author{Carlos H. Villanueva	\and
		Kurt Maute
}

\institute{ C. H. Villanueva \at
				Department of Mechanical Engineering,\\
				University of Colorado at Boulder,\\
				Boulder, CO 427 UCB, USA\\
				e-mail: carlos.villanueva@colorado.edu\\
			\and
			K. Maute \at
				Department of Aerospace Engineering,\\
				University of Colorado at Boulder,\\
				Boulder, CO 427 UCB, USA\\
				e-mail: maute@colorado.edu\\
}

\date{Received: date / Accepted: date}

\maketitle

%%%%%%%%%%%%%%%%%%%%%%%%%%%%%%%%%%%%%%%%%%%%%%%%%%%%%%%%%%%%%%%%%%%
% Abstract
%%%%%%%%%%%%%%%%%%%%%%%%%%%%%%%%%%%%%%%%%%%%%%%%%%%%%%%%%%%%%%%%%%%

\begin{abstract}
As the capabilities of additive manufacturing techniques increase, topology optimization provides a promising approach to design geometrically sophisticated structures which can be directly manufactured. Traditional topology optimization methods aim at finding the conceptual design but often lack a sufficient resolution of the geometry and structural response, needed to directly use the optimized design for manufacturing. To overcome these limitations, this paper studies the viability and characteristics of the eXtended Finite Element Method (XFEM) in combination with the Level-Set Method (LSM) for topology optimization of three dimensional structural design problems. The LSM describes the geometry by defining the nodal level set values via explicit functions of the optimization variables. The structural response is predicted by a generalized version of the XFEM. The LSM-XFEM approach is compared against results from a traditional Solid Isotropic Material with Penalization (SIMP) method for two-phase ``solid-void'' and ``solid-solid'' problems. The numerical results demonstrate that the LSM-XFEM approach can describe crisply the geometry and predict the structural response of complex three-dimensional structures with acceptable accuracy even on coarse meshes. However, the LSM-XFEM studied here lacks a robust and intuitive formulation to control the minimum feature size, and the optimization results may depend on the initial design.


\end{abstract}

%%%%%%%%%%%%%%%%%%%%%%%%%%%%%%%%%%%%%%%%%%%%%%%%%%%%%%%%%%%%%%%%%%%
% Introduction
%%%%%%%%%%%%%%%%%%%%%%%%%%%%%%%%%%%%%%%%%%%%%%%%%%%%%%%%%%%%%%%%%%%

\section{Introduction}
\label{sec:intro}

Recent advances in additive manufacturing allow the precise placement of one or multiple materials at micrometer resolution with essentially no restrictions on the geometric complexity of the spatial arrangement. Complex three dimensional solids can be created with highly non-regular material distributions in a near optimal fashion, enabling the fabrication of structures with enhanced performance. Topology optimization has emerged as a promising approach to utilize the benefits of additive manufacturing \citep{NP:12,MGW+:13}. Structural topology optimization seeks to find the optimal geometry and/or the material layout of a body within a given design domain. The geometry is represented by the spatial distribution of two or more material phases; in structural problems, one of these material phases may represent void.

Originally topology optimization methods were developed primarily to create conceptual designs in the early stage of the design process \citep{BS:03,Rozvany:09}. Later, topology optimization was applied to directly design micro-electro-mechanical systems (MEMS), utilizing the ability of thin-film fabrication techniques, such as photolithography in combination with chemical etching, to create geometrically complex devices at low cost \citep{Sigmund:01,Sigmund:01a}. As the MEMS design problem is essentially two-dimensional, traditional approaches were sufficient to achieve the necessary geometric resolution at acceptable computational cost. Motivated by the availability of affordable additive fabrication methods for three dimensional structures, this paper focuses on topology optimization of three dimensional structures and introduces a new approach for finding optimized designs with high geometric resolution on rather coarse computational meshes.      

Most approaches for structural topology optimization are density methods. For a two-phase problem, the density is considered a design variable and can assume intermediate values between the density of the material phase ``A'' and the density of the material phase ``B''. The most popular density method is the Solid Isotropic Material with Penalization method introduced by \citet{Bendsoe:89} and \citet{ZR:91}. It features great versatility, robustness, efficiency, and ease of use and implementation for a broad range of applications \citep{SM:13,DG:13}.

Density methods typically describe the boundaries between the material phases either via intermediate density values or by discrete material distributions leading to jagged boundary geometries. In both cases, the enforcement of boundary and interface conditions at the material interface is hampered and may result in non-physical responses, such as premature yielding \citep{MSR:98}. Often this issue can be mitigated by mesh refinement or adaptive re-meshing \citep{MR:95,MR:97}. However, for problems that require an accurate description of the boundaries, such as boundary layer problems in fluids and skin-depth issues in electromagnetics, it was reported that density methods like SIMP fail \citep{SM:13}.

The shortcomings of density methods have promoted the development of the Level Set Method (LSM) for topology optimization. The LSM allows a crisp representation of the phase boundaries and the accurate enforcement of boundary conditions on fixed meshes. The material interface in the LSM is described implicitly by the iso-contours of a Level Set Function (LSF), usually the zero level-set contour \citep{AJT:04,SW:00,WWG:03}.

The LSF is typically discretized by the same mesh used for the physical field and is updated in the optimization process via the solution of the Hamilton-Jacobi equations. Alternatively, the parameters of the discretized LSF are defined as explicit functions of the optimization variables, and the resulting parameter optimization problem is solved by standard nonlinear programming methods \citep{DML+:13}. The key challenges for the LSM include (a) controlling the spatial gradients of the LSF in the vicinity of the zero-level set contour to avoid ill-conditioning of the optimization problem, (b) controlling local feature sizes, and (c) accelerating the convergence of the geometry in the optimization process. For a detailed discussion of the LSM, the reader is referred to the comprehensive review by \citet{DML+:13} and \citet{GP:13}.

In LSMs, the structural geometry can be represented in the mechanical model via an Ersatz material approach, immersed boundary techniques, or by adaptive geometry conforming meshes. The first two approaches allow the use of fixed, design independent meshes while the last approach requires local or global re-meshing as the structural geometry evolves in the design process.

Using an Ersatz material approach, the void phase is modeled by a soft material and the material properties in elements intersected by the zero-level set contour are interpolated between the ``void'' and solid phase, proportional to the volume ratio of the individual phases. However, this approach faces the same issues as density methods in regards to enforcing boundary conditions across the material interface. Other approaches to model the mechanical response include generalized and adaptive finite element schemes such as the Super-Imposed Finite Element Method (SFEM) \citep{WW:06}, the eXtended finite element method (XFEM) \citep{MD:07,WWX:10,KM:12}, and local re-meshing schemes \citep{YNK+:11}.

In this paper, we focus on the LSM in combination with the XFEM. The XFEM does not require a mesh that conforms to the material interfaces and reduces the complexity of mesh construction. Spatial discontinuities in the structural response are captured by augmenting the standard finite element interpolations with additional shape functions. This approach is similar to the SFEM \citep{WW:06} as the solution is obtained by super-imposing the standard and enriched shape functions. However, unlike the SFEM, the XFEM can combine multiple types of shape functions and thus allows for greater flexibility.

The XFEM builds upon the partition of unity concept developed by \citet{NME:NME86}. The idea of XFEM was originally proposed by \citet{BB:99} to model crack propagation. The reader is referred to \citet{ANH:09} for an overview of the application of XFEM to problems in fracture mechanics. The XFEM has been used for a variety of problems in computational mechanics, such as fluid-structure interaction \citep{GW:08,GW:08a}, multi-phase flows \citep{Fries:09}, and nano-scale heat transfer \citep{LYM:11}. A general overview of the method is presented by \citet{FB:10}.

\citet{DMJ+:06} originally introduced the XFEM into shape optimization using a simplified XFEM formulation. This formulation does not use additional enrichment functions and is limited to problems where one of the material phases represents void and the geometric configuration is ``simple'', i.e.~does not contain geometric features that are smaller than the size of two elements \citep{MM:13}. In this instance, the weak form of the governing equations is only integrated over the solid material in each element. In addition, if the interface between the material phases is traction free, this simplified version of the XFEM only differs from the traditional finite element method with respect to the domain of integration. The simplified XFEM version was applied to structural shape optimization, for example, by \citet{DMJ+:06}, \citet{MMF+:05}, and \citet{MD:07}, and to topology optimization of three dimensional structures by \citet{LWW:12}.

An XFEM approach based on a standard enrichment strategy allows to model two-phase problems with a simple intersection pattern. This formulation is considered, for example, by \citet{WWX:10} to solve ``solid-void'' structural topology optimization problems, modeling the ``void'' as a soft material. \citet{MKM+:11} used a standard enrichment strategy to discretize the phonon Boltzmann transport equation and optimize the thermal conductivity of nano-composites. However, this enrichment strategy is not guaranteed to consistently approximate the state variable fields for configurations with complex intersection patterns.

An enhanced version of the XFEM was proposed by \citet{MM:13}, who presented a generalized enrichment based on the step enrichment of \citet{HH:04} and applied it to two dimensional structural topology optimization. This formulation captures consistently the mechanical response for complex geometries and intersection patterns for general multi-phase problems.

This paper will expand the combination of the LSM and the XFEM onto general two-phase, three dimensional problems. We will compare results of the proposed LSM-XFEM framework with SIMP results for structural topology optimization examples. The numerical results will show that the LSM-XFEM combination is a promising approach for three dimensional problems and allows for the use of coarse meshes to represent the structural geometry and to describe the structural response with acceptable accuracy.

The main challenges of expanding the previous LSM with XFEM approaches to three dimensions stem from the increased complexity of possible intersections patterns. Such patterns include elements that are intersected multiple times and elements containing only a small volume of a particular phase. In contrast to \citet{LWW:12}, we will adopt the generalized enrichment strategy of \citet{MM:13} to (a) accurately model the structural response on complex three dimensional patterns and (b) solve solid-solid material distribution problems. Further, we will expand the preconditioning scheme of \citet{LMD+:13} onto three dimensions to mitigate ill-conditioning issues in the XFEM analysis problems due to elements with small volume fractions.

The main contributions of the paper are: (a) we present a numerically robust and computationally viable approach for solving general two-phase, three dimensional topology optimization problems, and (b) we provide a direct comparison of LSM-XFEM and SIMP results for three dimensional problems, highlighting key features of the two methods.

The remainder of this paper is structured as follows: the geometry models of the LSM-XFEM and SIMP methods are\ described in Section \ref{sec:geometry-modeling}. The mechanical model and the XFEM formulation are summarized  in Section \ref{sec:structural-analysis}. Details of the LSM-XFEM and SIMP optimization approaches are presented in Section \ref{sec:optimization-model}. Section \ref{sec:computational-considerations} highlights specific computational challenges of the XFEM approach for three dimensional problems. Section \ref{sec:numerical-examples} presents structural topology optimization examples, comparing the LSM-XFEM and SIMP approaches. Section \ref{sec:conclusions} summarizes the main conclusions drawn from this study.

%%%%%%%%%%%%%%%%%%%%%%%%%%%%%%%%%%%%%%%%%%%%%%%%%%%%%%%%%%%%%%%%%%%
% Geometry Modeling
%%%%%%%%%%%%%%%%%%%%%%%%%%%%%%%%%%%%%%%%%%%%%%%%%%%%%%%%%%%%%%%%%%%

\section{Geometry Modeling}
\label{sec:geometry-modeling}

In topology optimization, the geometry of a body is defined via its material distribution. In density methods, such as SIMP, the material distribution is discretized by finite elements, with either elemental or nodal parameters defining the distribution within the element. Most often the same mesh is used to approximate the density distribution and the structural response. Alternatively, the state and density fields can be discretized by different meshes with different refinement levels; see for example the Multi-resolution Topology Method (MTOP) by \citet{NSP:10}. The optimization variables define analytically or by means of auxiliary partial differential equations the nodal or element density parameters \citep{SM:13}. For two-phase problems, the density is continuously varied between ``0'' (phase ``A'') and ``1'' (phase ``B''). Implicit or explicit penalization schemes, optionally combined with projection methods, are used to encourage ``0-1'' solutions \citep{GPB:04,Sigmund:07}.

The crispness of the interface geometry, as described via the optimized material distribution, depends on (a) the formulation of the optimization problem, i.e.~the objective and constraints, (b) regularization techniques, such as density or sensitivity filters, and (c) the optimization algorithm. In general, the resolution of the phase boundaries increases as the mesh is refined. For three dimensional problems, often coarse meshes are used to limit the computational costs. In this case, the boundary geometry either lacks crispness due to the presence of elements with intermediate densities or is approximated by a spatially discontinuous material distribution, leading to jagged interfaces.

Alternatively, the material distribution can be described via a level set function $\phi(\mathbf{x})$. Typically the zero level set contour defines the phase boundaries. Considering a two-phase problem, the interface $\Gamma_{A,B}$ is defined as follows:

\begin{equation}
\label{eq:level-set}
	\begin{aligned}
		\phi(\mathbf{x}) &< 0,    && \forall \ \mathbf{x} \in \Omega_A, \\
		\phi(\mathbf{x}) &> 0,    && \forall \ \mathbf{x} \in \Omega_B, \\
		\phi(\mathbf{x}) &= 0,    && \forall \ \mathbf{x} \in \Gamma_{A,B},
	\end{aligned}
\end{equation}

where the vector $\mathbf{x}$ collects the spatial coordinates, $\Omega_{A}$ is the domain of material phase ``A'', $\Omega_{B}$ is the domain of material phase ``B'', and $\Gamma_{A,B}$ defines the material interface between phase ``A'' and phase ``B''. For example, to model a sphere in a three dimensional mesh centered at ($x_c$,$y_c$,$z_c$), the value of the level set function at a grid point ($x_i$,$y_i$,$z_i$) is:

\begin{equation}
	\label{eqn:circle}
	\phi_i = { (x_i - x_c) }^2 + { (y_i - y_c) }^2 + { (z_i - z_c) }^2- r^2,
\end{equation}

where $r$ represents the radius of the sphere, and the sign value of $\phi_{i}$ at each node determines if the node is inside or outside the circle.

The level set field is typically discretized by shape functions with either local or global support; the reader is referred to the review paper by \citet{DML+:13}. In the simplest and most common approach, the level-set field is approximated on the same mesh used for discretizing the governing equations. In this study, we follow this approach and define the nodal level-set values explicitly as functions of the optimization variables (see \ref{sec:optimization-XFEM}). The resulting parameter optimization problem is solved by a standard nonlinear programming method.

\begin{figure*}
	\includegraphics[width=\linewidth]{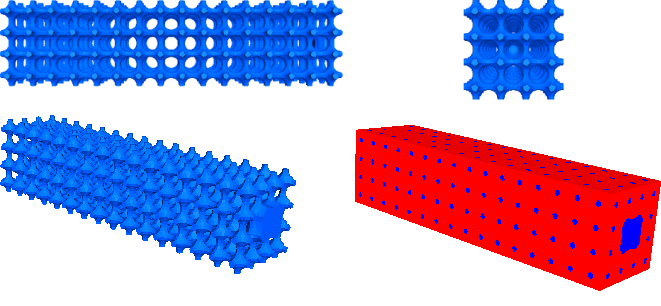}
	\caption{Initial design with array of spherical inclusions for the cantilever beam example of Section \ref{sec:two-phase-optimization}.}
	\label{fig:long-cantilever-beam-swiss-cheese}
\end{figure*}

While LSMs provide a crisp definition of the phase boundaries, they require seeding the initial design with inclusions and/or introducing inclusions in the course of the optimization process, for example via topological derivatives \citep{EKS:94,SZ:99,BHR:04,NBH+:07}. An example of an initial design with a regular pattern of spherical inclusions is shown in Fig.~\ref{fig:long-cantilever-beam-swiss-cheese}. The images in the upper row and the image in the lower left corner show only phase ``A''. The material layout of both phases is depicted in the lower right image. This layout can be generated by superposing Eq. \ref{eqn:circle} for spheres at uniformly spaced center locations.

The optimization results of the LSM are typically dependent on the initial layout. Furthermore, it is non-trivial to generate an initial design that satisfies geometric design constraints, such as volume or perimeter constraints. In our experience, severe constraint violations may cause the optimization process to converge to a feasible design with otherwise poor performance.

%%%%%%%%%%%%%%%%%%%%%%%%%%%%%%%%%%%%%%%%%%%%%%%%%%%%%%%%%%%%%%%%%%%
% Structural Analysis
%%%%%%%%%%%%%%%%%%%%%%%%%%%%%%%%%%%%%%%%%%%%%%%%%%%%%%%%%%%%%%%%%%%

\section{Structural Analysis}
\label{sec:structural-analysis}

In this section, we briefly discuss the structural model and the XFEM analysis used in this paper. The governing equations are presented first, followed by a summary of the XFEM discretization and analysis.

%%%%%%%%%%%%%%%%%%%%%%%%%%%%%%%%%
% Structural Model
%%%%%%%%%%%%%%%%%%%%%%%%%%%%%%%%%

\subsection{Structural Model}
\label{sec:structural-model}

\begin{figure}
	\includegraphics[width=0.8\linewidth]{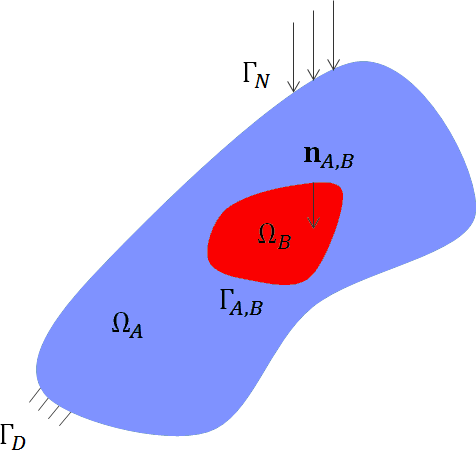}
	\caption{XFEM two-phase problem model.}
	\label{fig:xfem-problem-setup}
\end{figure}

This paper considers the topology optimization of structures using the LSM approach described above and the XFEM to predict the structural response, assuming infinitesimal strains, a linear elastic material behavior, and static conditions. We consider the two-phase problem depicted in Fig.~\ref{fig:xfem-problem-setup}, where $\Gamma_{N}$ denotes the surface where traction forces are applied, $\Gamma_{D}$ denotes the surface with prescribed displacements, and $\mathbf{n}$ is the normal at the material interface pointing from phase ``A'' to phase ``B''. The weak form of the governing equations can be decomposed into the following terms:

\begin{equation}
	\label{eq:residual-equation}
	W = W_{S} + W_{L} + W_{k} = 0,
\end{equation}

where $W_{S}$ collects the contributions from the static equilibrium, including body forces and surface tractions, $W_{L}$ models the interface conditions along the phase boundaries for ``solid-solid'' problems, and $W_{k}$ is due to a fictitious spring model to pin free floating material in ``solid-void problems''.

The weak form of the structural equilibrium equations is:

\begin{multline}
\label{eq:weak-form}
	W_{S} = \int_{\Omega_A} { \pmb{\eta} \colon \pmb{\sigma}(\mathbf{u}) \ d \Omega}
		  + \int_{\Omega_B} { \pmb{\eta} \colon \pmb{\sigma}(\mathbf{u}) \ d \Omega}
		  - \int_{\Omega_A} { \mathbf{v} \cdot \mathbf{b} \ d \Omega } \\
		  - \int_{\Omega_B} { \mathbf{v} \cdot \mathbf{b} \ d \Omega }
		  - \int_{\Gamma_N} { \mathbf{v} \cdot \mathbf{f} \ d \Gamma_N },
\end{multline}

where $\mathbf{v}$ is the kinematically admissible test function, $\pmb{\eta}$ is the strain tensor associated with the test function $\mathbf{v}$, $\mathbf{u}$ is the displacement vector, $\pmb{\sigma}$ is the stress tensor, $\mathbf{b}$ is the applied body force, and $\mathbf{f}$ is the external traction applied along $\Gamma_{N}$.

To enforce continuity of the displacements along the phase boundaries, the static equilibrium equations are typically augmented by either an enhanced Lagrange multiplier or penalty formulations, such as the stabilized Lagrange multiplier and the Nitsche method. Note that the standard Lagrange multiplier approach is not suitable for the XFEM as it suffers from numerical instabilities. The reader is referred to \citet{Stenberg:95}, \citet{JS:09}, and \citet{DH:09} for more details.

In this paper, we enforce displacement continuity along phase boundaries for ``solid-solid'' problems via the following stabilized Lagrange multiplier method:

\begin{multline}
\label{eq:stabilized_Residual}
	W_{L} = - \int_{\Gamma_{A,B}} \left[ \mathbf{v} \right] \cdot \pmb{\lambda} \ d \Gamma_{A,B}
		  + \gamma \int_{\Gamma_{A,B}} \pmb{\mu} \cdot \left[ \mathbf{u} \right] d \Gamma_{A,B} \\
		  + \int_{\Gamma_{A,B}} \mu \cdot \left( \pmb{\lambda} - \bar{\pmb{\sigma}} \cdot \mathbf{n}_{A,B} \right) d \Gamma_{A,B},
\end{multline}

\begin{equation} \label{jump}
	\left[ \mathbf{u} \right] = \mathbf{u}^{(A)} - \mathbf{u}^{(B)} \ ,
	\left[ \mathbf{v} \right] = \mathbf{v}^{(A)} - \mathbf{v}^{(B)} \ ,
\end{equation}

\begin{equation}
\label{eq:lagrangeStress}
	\bar{\pmb{\sigma}} = \frac{1}{2} \left( \pmb{\sigma}^{(A)} + \pmb{\sigma}^{(B)} \right),
\end{equation}

where $\pmb{\lambda}$ is the Lagrange multiplier, and $\pmb{\mu}$ is the associated test function. The higher the weight $\gamma$ is, the better the interface condition is satisfied, at the cost of numerical stability.

The formation of free floating solid particles surrounded by void material is possible in ``solid-void'' topology optimization, leading to a singular analysis problem. This issue does not exist in an Ersatz material approach as the void phase is modeled via a soft material. A similar approach can be applied to the XFEM to suppress singularities by modeling the ``void'' phase via a soft material \citep{WWX:10}. However, this approach requires accounting for the interface contributions (\ref{eq:stabilized_Residual}) and integrating the governing equations over the void phase. To avoid the associated complexity and computational costs, we extend the approach of \citet{MM:13} onto three dimensions and assume that the solid phase is supported by fictitious springs. This model leads to the following contribution to the governing equations, assuming that phase ``A'' is the solid phase:

\begin{equation}
\label{eq:residual-k}
	W_{k} = \int_{\Omega_A} { k \ \mathbf{v} \cdot \mathbf{u} \ d \Omega },
\end{equation}

where $k$ denotes the stiffness of the distributed system of springs.

For a more detailed explanation of this XFEM formulation, the reader is referred to the paper by \citet{MM:13}.

%%%%%%%%%%%%%%%%%%%%%%%%%%%%%%%%%
% Discretization
%%%%%%%%%%%%%%%%%%%%%%%%%%%%%%%%%

\subsection{Discretization}
\label{sec:discretization}

To capture the discontinuities in the strain and stress fields along the phase boundaries, we enrich the standard finite element approximation with additional shape functions. We adopt the generalized enrichment strategy of \citet{MM:13} which resolves consistently the displacement fields in the presence of small features and does not suffer from artificially coupling disconnected phases. Considering a two-phase problem, the displacement field is approximated as follows:

\begin{equation}
\label{eq:xfem-interp}
	\mathbf{u}(\mathbf{x}) = \sum \limits^{M}_{m=1} \left( H(-\phi) \sum\limits^{n}_{i=1} \mathbf{N}_i \ \mathbf{u}_{i,m}^A
														 + H( \phi) \sum\limits^{n}_{i=1} \mathbf{N}_i \ \mathbf{u}_{i,m}^B \right),
\end{equation}

where $m$ is the enrichment level, $M$ is the maximum number of enrichment levels used for each phase, $\mathbf{N}$ are the shape functions, $\mathbf{u}^l_{i,m}$ is the vector of nodal displacement components at node $i$ for phase $l=[A,B]$, $\phi$ is the level set value evaluated at the integration point, and $H$ denotes the Heaviside function. The enrichment level is chosen such that the displacements in disconnected volumes of the same phase are interpolated by separate sets of degrees of freedom. When interpolating the level set field by element-wise linear functions, a maximum of $14$ enrichment levels is needed in three dimensions. This enrichment strategy will be revisited in Section \ref{sec:computational-considerations}.

The Heaviside function $H$ depends on the level set function and is defined as follows:

\begin{equation}
\label{eq:heaviside}
	H(z) =
		\begin{cases}
			1 & z > 0, \\
			0 & z \le 0,
		\end{cases}
\end{equation}

The Heaviside functions ``turns on/off'' the standard finite element interpolations in the particular phases. The approximation (\ref{eq:xfem-interp}) allows for discontinuities of the displacements along the phase boundaries. Therefore the continuity is enforced weakly via the stabilized Lagrange multiplier method (\ref{eq:stabilized_Residual}).

Following a Bubnov-Galerkin scheme, we test the governing equations with the same subspace as we use for the trail functions; see Eq. \ref{eq:xfem-interp}. The weak form of the governing equations is integrated numerically over the individual phases, using the Delaunay triangulation of the element along the phase boundaries.

%%%%%%%%%%%%%%%%%%%%%%%%%%%%%%%%%
% Preconditioner
%%%%%%%%%%%%%%%%%%%%%%%%%%%%%%%%%

\subsection{Preconditioner}
\label{sec:preconditioner}

As described above, the degrees of freedom $\mathbf{u}^l_{i,m}$ interpolate the structural displacements in topologically connected subdomains of phase $l$ in the elements connected to node $i$. As the total volume of these subdomains vanishes, the discretized structural model becomes increasingly ill-conditioned; i.e.~the condition number of the stiffness matrix rapidly increases. This phenomenon is more pronounced for three dimensional problems than those in two dimensional ones.

To mitigate this ill-conditioning issue, we extend the geometric preconditioning scheme of \citet{LMD+:13}, which was introduced and studied for two dimensional heat conduction and flow problems, onto three dimensional problems in structural mechanics. The goal of this preconditioning scheme is to balance the influence of all degrees of freedom in the system, as the volumes in which the subset of these degrees of freedoms interpolates the solution approach zero. To this end, we introduce the following projection:

\begin{equation}
\label{eqn:preconditioner-u}
	\tilde{\mathbf{u}} = \pmb{T} \mathbf{u},
\end{equation}

where $\mathbf{u}$ is the vector of displacement degrees of freedom according to Eq. \ref{eq:xfem-interp}, $\pmb{T}$ is a transformation matrix,  and $\tilde{\mathbf{u}}$ is the solution vector in the transformed space. The residual, $\tilde{\pmb{R}}$, and stiffness matrix, $\tilde{\pmb{K}}$, in the transformed space are defined as:

\begin{equation}
\label{eqn:preconditioner-residual}
	\tilde{\pmb{R}} = \pmb{T}^{T} \pmb{R}  \qquad 	\tilde{\pmb{K}} = \pmb{T}^{T} \pmb{K} \pmb{T},
\end{equation}

where the residual, $\pmb{R}$, and the stiffness matrix, $\pmb{K}$, result from integrating the weak form of the governing equations using the XFEM approximation (\ref{eq:xfem-interp}).

The preconditioner $\pmb{T}$ is a diagonal matrix built by integrating the spatial derivatives of the shape functions over the nodal support of nodes connected to an intersected element. The diagonal components of the matrix are defined as

\begin{equation}
\label{eqn:T-matrix}
	\pmb{T}_{i,m}^{l} = \left( \max_{e \in E_{i}} \frac { \int_{\mathcal{D}_l^e} \nabla \mathbf{N}_i(\mathbf{x}) \cdot
	\nabla
	\mathbf{N}_i(\mathbf{x}) \,d\mathbf{x}}{\int_{\mathcal{D}^e} \nabla \mathbf{N}_i(\mathbf{x}) \cdot \nabla
	\mathbf{N}_i(\mathbf{x}) \,d\mathbf{x}} \right)^{-1},
\end{equation}

where $\pmb{T}_{i,m}^{l}$ corresponds to the degree of freedom $\mathbf{u}_{i,m}^{l}$, $i$ is the node index, $l$ is the material phase ``A'' or ``B'', $m$ is the enrichment level, $E_{i}$ is the set of elements connected to node $i$, and $\mathcal{D}_l^e$ is the element domain of phase $l$. The components of the matrix increase as the region of influence of a degree of freedom decreases. The entries $T_{i,m}^{l}$ of nodes $i$ that are not connected to at least one intersected element are set to one.

To avoid numerical issues due to large values for the components of $\pmb{T}$, the degrees of freedom associated with the diagonal entry $\pmb{T}_{i,m}^{l}$ are constrained to zero if the following condition is satisfied:

\begin{equation}
\label{eqn:T-matrix-constraint}
	\pmb{T}_{i,m}^{l} \ge T_{tol},
\end{equation}

where $T_{tol}$ is a specified tolerance. As studies by \citet{LMD+:13} have shown, the above preconditioning scheme is rather insensitive to the value of $T_{tol}$ and is typically set to a value larger than $10^8$. For more details about the specifics of the formulation, the reader is referred to the paper by \citet{LMD+:13}.

%%%%%%%%%%%%%%%%%%%%%%%%%%%%%%%%%%%%%%%%%%%%%%%%%%%%%%%%%%%%%%%%%%%
% Optimization Model
%%%%%%%%%%%%%%%%%%%%%%%%%%%%%%%%%%%%%%%%%%%%%%%%%%%%%%%%%%%%%%%%%%%

\section{Optimization Model}
\label{sec:optimization-model}

The design optimization problems considered in this paper can be written as follows:

\begin{equation}
\label{eq:topo_problem}
	\begin{aligned}
		&\min_{\mathbf{s}} \mathcal{F}(\mathbf{s}, \mathbf{u}(\mathbf{s})), \\
		&\text{s.t.\,}
		\begin{cases}
			\mathbf{s}, & \text{subject to design constraints}  \ \mathcal{G}_j \leq 0 \text{,}\\
			\mathbf{u}, & \text{solves} \ W= 0 \ \text{for a given $\mathbf{s}$,}
		\end{cases}
	\end{aligned}
\end{equation}

where $\mathbf{s}$ denotes the vector of design variables, $\mathcal{F}$ the objective function, and $\mathcal{G}_j$ the $j$-th design constraint. In general, the objective and constraints depend on the optimization and state variables. The optimization problem (\ref{eq:topo_problem}) is solved by nonlinear programming methods, and the gradients of the objective and constraint functions are computed via the adjoint method.

In this paper, we compare the proposed LSM-XFEM against the well-known SIMP method, augmented by a projection scheme. In the following subsection we briefly outline the models that define the discretized level set field (LSM) and the material properties (SIMP) as function of the optimization variables.

%%%%%%%%%%%%%%%%%%%%%%%%%%%%%%%%%
% XFEM
%%%%%%%%%%%%%%%%%%%%%%%%%%%%%%%%%

\subsection{XFEM}
\label{sec:optimization-XFEM}

The nodal values of the discretized level set field are defined as analytical functions of the optimization variables via the following linear filter:

\begin{equation}
\label{eq:level-set-filter}
	\phi^{n}(\mathbf{s}) = \frac{ \sum \limits^{P}_{i=1} \mathbf{w}^{n}_{i} \mathbf{s}_i } { \sum \limits^{P}_{i=1} \mathbf{w}^{n}_{i} },
\end{equation}

\noindent
with

\begin{equation}
\label{eq:level-set-filter-weight}
	\mathbf{w}^{n}_{i} = max( 0, \parallel \mathbf{x}_i - \mathbf{x}_n \parallel - r_{\phi} ),
\end{equation}

where $P$ is the number of nodes in the discrete model, $\mathbf{x}_{i}$ is the location of the node at which the design variable $i$ is defined, $r_{\phi}$ is the filter radius, $\mathbf{w}^{n}_{i}$ is the weight of node $n$ with respect to design variable $i$, and $\mathbf{x}_n$ and $\phi^{n}$ are the position vector and the computed level set value of node $n$, respectively.

The above linear filter was used previously in the studies of \citet{KM:12} and \citet{MM:13}, and was shown to improve the convergence rate in the optimization process. However, in contrast to density or sensitivity filters used in SIMP methods, the filter above is not guaranteed to control the minimum feature size. This issue will be revisited in Section \ref{sec:numerical-examples}.

%%%%%%%%%%%%%%%%%%%%%%%%%%%%%%%%%
% SIMP
%%%%%%%%%%%%%%%%%%%%%%%%%%%%%%%%%

\subsection{SIMP}
\label{sec:optimization-SIMP}

Here, the material distribution is parameterized by nodal density values, $\rho_i$, which are treated as optimization variables, i.e.~$\rho_i=\mathbf{s}_i$. Following the work of \citet{GPB:04}, we compute the elemental density by combining a linear density filter and a projection scheme as follows:

\begin{equation}
\label{eq:smoothing-equation-SIMP}
	\rho^{e}(\mathbf{s}) = \frac{ \sum \limits^{E}_{i=1} \mathbf{w}^{e}_{i} \mathbf{s}_i }
					{ \sum \limits^{E}_{i=1} \mathbf{w}^{e}_{i}     },
\end{equation}

\noindent
where

\begin{equation}
\label{eq:smoothing-radius-SIMP}
	\mathbf{w}^{e}_{i} = max( 0, \parallel \mathbf{x}_i - \mathbf{x}_e \parallel - r_{\rho} ),
\end{equation}

where $E$ is the number of elements in the discrete model, $\mathbf{x}_{i}$ is the location of the node at which the design variable $i$ is defined, $r_{\rho}$ is the filter radius, $\mathbf{w}^{e}_{i}$ is the factor of element $e$ with respect to design variable $i$, and $\mathbf{x}_e$ and $\rho^{e}$ are the position vector of the centroid and the computed elemental density of element $e$, respectively.

\citep{GPB:04} proposed a density projection method to reduce the volume occupied by material with intermediate densities. The projection is based on a smoothed Heaviside function and applied to the elemental densities as follows:

\begin{equation}
\label{eq:heaviside-SIMP}
	\hat{\rho}^{e}(\mathbf{s}) = 1 - e^{-\beta \rho^{e}(\mathbf{s})} + \rho^{e}(\mathbf{s}) e^{-\beta}
\end{equation}

where $\hat{\rho}^{e}$ is the projected elemental density, and the parameter $\beta \ge 0$ controls the crispness of the projection. For $\beta = 0$ the projection turns into an identity operator, i.e.~$\hat{\rho}^{e}= \rho^{e}$.

The Young's modulus, $E$, is defined as a function of the density, $\hat{\rho}^{e}$, using the standard SIMP interpolation:

\begin{equation}
\label{eq:SIMP-elastic-modulus-two-phase}
	E(\mathbf{x}) = E_{B} + \lvert E_{A} - E_{B} \rvert \hat{\rho}^{e}(\mathbf{s})^{p}
\end{equation}

where $E_{A}$ and $E_{B}$ are the Young's moduli for material phase ``A'' and ``B'', and $p$ is the SIMP penalization factor. To model a ``solid-void'' optimization problem, $E_B$ is set to value much smaller than $E_A$.

The filter (\ref{eq:smoothing-equation-SIMP}) prevents the formation of features smaller than $r_{\rho}$, typically at the cost of generating intermediate density values along the phase boundaries. This effect is mitigated by the projection (\ref{eq:heaviside-SIMP}) which, in the limit for $\beta \rightarrow \infty$ , maps non-zero $\rho^e$ values into ``1''. The reader is referred to the papers by \citet{GPB:04} and \citet{GAH:11} for further details of the scheme presented above, and to \citet{SM:13} for a comprehensive discussion of projections schemes.

%%%%%%%%%%%%%%%%%%%%%%%%%%%%%%%%%%%%%%%%%%%%%%%%%%%%%%%%%%%%%%%%%%%
% Computational Considerations
%%%%%%%%%%%%%%%%%%%%%%%%%%%%%%%%%%%%%%%%%%%%%%%%%%%%%%%%%%%%%%%%%%%

\section{Computational Considerations}
\label{sec:computational-considerations}

\begin{figure}
	\includegraphics[width=\linewidth]{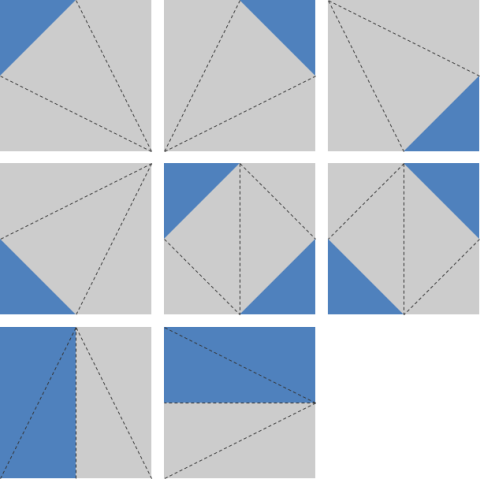}
	\caption{Intersection pattern for a two-dimensional QUAD4 element.}
	\label{fig:2D-intersections}
\end{figure}

\begin{figure}[hb!]
	\includegraphics[width=\linewidth]{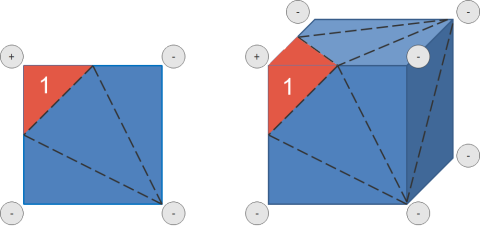}
	\caption{Initial enrichment function computation.}
	\label{fig:subphase-2D-to-3D}
\end{figure}

Expanding the LSM-XFEM combination onto three dimensional problems faces both algorithmic and computational challenges which are briefly discussed below.

\begin{figure*}[t]
	\includegraphics[width=\linewidth]{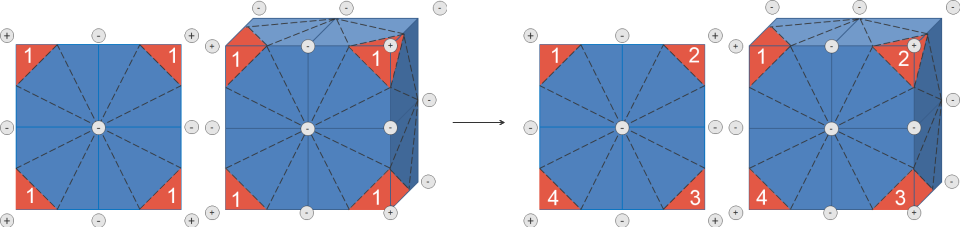}
	\caption{Distribution of the enrichment functions after the application of the enrichment strategy.}
	\label{fig:enrichment-2D-to-3D}
\end{figure*}

The XFEM requires integrating the weak form of the governing equations separately in each phase. To this end, an element intersected by the zero level set contour is subdivided. For two dimensional problems and using a linear interpolation of the level set field within an element, there are only 8 intersection configurations which can be tabulated; see Fig.~\ref{fig:2D-intersections}. In three dimensions, there are 127 intersection configurations. To handle this complexity, we compute the intersection point of the zero level set contour with the element edges and use a Delaunay triangulation to subdivide the element. In numerical experiments, this approach has proven robust and computationally inexpensive.

Previous studies on topology optimization for three dimensional structures with the XFEM \citep{LWW:12,HMM:13} have employed a simplified enrichment scheme which is limited to ``solid-void'' problems and may suffer from artificial coupling of disconnected material. Our work overcomes these issues by adopting the generalized enrichment scheme summarized in Section \ref{sec:discretization}. The key challenge of this scheme is to identify the enrichment levels needed to consistently interpolate the displacements in elemental subdomains with the same phase.

To this end, the subdomains in all elements connected to a node need to be considered. This naturally leads to an algorithm which loops over all nodes and, in an inner loop, over all elements connected to the current node. As this approach processes an element repeatedly, the following simple and efficient two-step scheme is introduced:

\begin{enumerate}

\item A temporary, elemental enrichment level is assigned to the subdomains in each element. Recall that the enrichment	level defines the set of degrees of freedom used to interpolate the displacements in an elemental subdomain. Note	that because this assignment is done individually for each element, the continuity of the interpolation across elements is not guaranteed. Fig.~\ref{fig:subphase-2D-to-3D} shows the triangulation and enrichment level for the red phase in two and three dimensions.

\item The nodal enrichment levels are constructed to ensure that the displacement field is interpolated continuously across elements, and by a different set of shape functions for each disconnected elemental subdomain of the same phase.  To this end, the cluster of elements connected to a node is considered, and the elemental enrichment levels assigned in step 1 are adjusted to satisfy the continuity and consistency conditions.
\end{enumerate}

This process is illustrated in Fig.~\ref{fig:enrichment-2D-to-3D}. The node of interest is the one located in the center of the element cluster. In step 1 each subdomain of the red phase is assigned an enrichment level of $m=1$. Applying this enrichment level to the degrees of freedom for the red phase at the center node would incorrectly couple the displacement fields in the red phase subdomains. Analyzing the element cluster around the center nodes shows that these subdomains are disconnected and individual enrichment level are assigned.

Topology optimization in three dimensions leads to FEM or XFEM models with a large number of degrees of freedom, and typically requires using iterative solvers and parallel computing. The stabilized Lagrange multiplier formulation of the interface conditions (\ref{eq:stabilized_Residual}) leads to a non-symmetric stiffness matrix. Numerical experiments have shown that the XFEM problems considered in this study can be robustly and efficiently solved by a generalized minimal residual (GMRES) method preconditioned by incomplete LU (ILU) factorization. Note that the ILU preconditioner operates on the projected XFEM system (\ref{eqn:preconditioner-residual}).

%%%%%%%%%%%%%%%%%%%%%%%%%%%%%%%%%%%%%%%%%%%%%%%%%%%%%%%%%%%%%%%%%%%
% Numerical Examples
%%%%%%%%%%%%%%%%%%%%%%%%%%%%%%%%%%%%%%%%%%%%%%%%%%%%%%%%%%%%%%%%%%%

\section{Numerical Examples}
\label{sec:numerical-examples}

We study the features of the proposed LSM-XFEM topology optimization approach with numerical examples. The LSM-XFEM results of ``solid-void'' and ``solid-solid'' problems are compared against the ones of the SIMP approach outlined in Section \ref{sec:optimization-SIMP}. In all examples we seek to minimize the strain energy subject to a constraint on the volume of the stiff phase. This problem formulation is chosen because it is well studied in the literature and the numerical experiments can be easily repeated. The following numerical studies will provide insight into (a) the convergence of the geometry and the structural response as the meshes are refined and (b) the influence of regularization techniques on the optimized results, such as the filter radii in (\ref{eq:level-set-filter}) and (\ref{eq:smoothing-equation-SIMP}), and perimeter constraints.

In all examples, the optimization problems are solved by the Globally Convergent Method of Moving Asymptotes (GCMMA) of \citet{Svanberg:02}. The sensitivities are computed by the adjoint method. The design domains are discretized by 8-node linear elements. The linear systems of the forward and adjoint problems are solved by a parallel implementation of the GMRES method \citep{Trilinos:03}. The problems are preconditioned by an ILU factorization with a fill of $2.0$ and an overlap of $1.0$. The convergence tolerances for both, the GCMMA and the GMRES solver, are chosen sufficiently low such that the optimization results do not depend on the tolerance values. In the LSM-XFEM examples, the spring stiffness value, $k$, is $10^{-6}$.

While the LSF-XFEM results can be directly used to fabricate the structure, for example by 3-D printing, the SIMP results need to be post-processed. From a practitioner perspective, only the post-processed SIMP results should be compared against the LSF-XFEM results. To this end we post-process the SIMP results with a lumping method that uses the iso-contour of the  density distribution. To obtain a ``0-1'' density distribution with smooth phase boundaries, we construct iso-contours for $\rho=\rho_T$ from the nodal density values, $\rho_i$; see Section \ref{sec:optimization-SIMP}. The volume enclosed by the iso-contour with $\rho \ge \rho_T$ is considered solid; the remaining volume is considered ``void''. The threshold, $\rho_T$, is determined such that the volume of the solid domain satisfies the volume constraint. The structural response of the results post-processed with the iso-contour density lumping (IDL) approach is analyzed conveniently with the XFEM.

To gain further insight into the crispness of the SIMP results and the influence of the post-processing methods above on their performance, we measure the volume fraction, $\bar{\rho}$, occupied by elemental densities with $0 < \hat{\rho}^e < 1$ as follows:

\begin{equation}
\label{eq:rho-utilization}
	\bar{\rho} = \frac{1}{ \int_{\Omega_D} {} \,d\Omega_D } \ \int_{\Omega_D} { \hat{\rho}^e  ( 1 - \hat{\rho}^e ) } \,
	d\Omega_D
\end{equation}

where $\Omega_D$ denotes the design domain.

%%%%%%%%%%%%%%%%%%%%%%%%%%%%%%%%%
% Cube with center load
%%%%%%%%%%%%%%%%%%%%%%%%%%%%%%%%%

\subsection{Cube with center load}
\label{sec:cube-center-load}

\begin{figure}[ht!]
	\centering
	\includegraphics[width=0.75\linewidth]{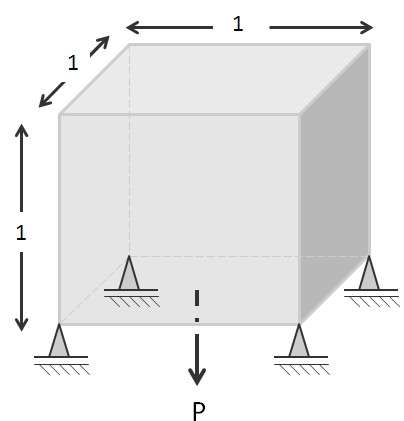}
	\caption{Cube with center load model.}
	\label{fig:cube-center-load-initial-setup}
\end{figure}

We consider the ``solid-void'' optimization problem depicted in Fig.~\ref{fig:cube-center-load-initial-setup}. With this example we will illustrate the basic features of the LSM-XFEM approach for three dimensional problems and show that the proposed LSM-XFEM approach and the SIMP formulation may exhibit comparable convergence behaviors as the mesh is refined.

The $1 \times 1 \times 1$ cubical design domain is pinned at its four bottom corners in the vertical direction and a unit force is applied at the center of the bottom face. The Young's modulus of the stiff phase is set to $1$ and the Poisson ratio to $0.3$. The maximum volume of the stiff phase is $10\%$. We compare LSM-XFEM and SIMP results for two mesh sizes: $24 \times 24 \times 24$ and $65 \times 65 \times 65$. The problem is solved on the full mesh.

First, we apply the SIMP approach with a penalization factor of $p=3$. The size of the smoothing radius is mesh dependent, and is set to $r_\rho = 3.2$ for the coarse mesh and $r_\rho = 1.182$ for the fine mesh; the projection parameter is set to $\beta=0$. Note that the smoothing radius is intentionally set relative to the element size ($1.6 \times$ the element edge). The Young's modulus of the void phase is set to $E_B = 10^{-9}$. While this approach does not ensure mesh-independent optimization results, it still prevents the formation of checker-board patterns and provides insight into the dependency of the geometry resolution of SIMP as the mesh is refined.

\begin{figure*}
	\begin{tabularx}{\linewidth}{XX}
		\subfloat[24x24x24 mesh]{\includegraphics[width=\linewidth]{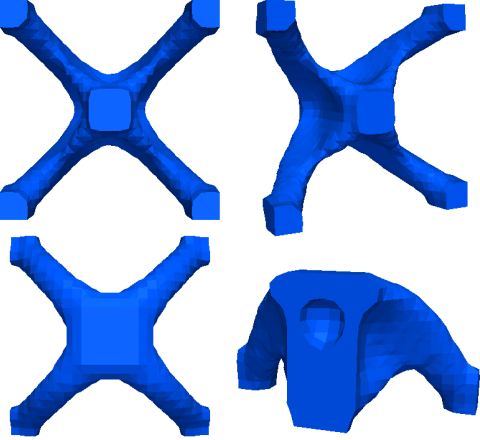}} &
		\subfloat[65x65x65 mesh]{\includegraphics[width=\linewidth]{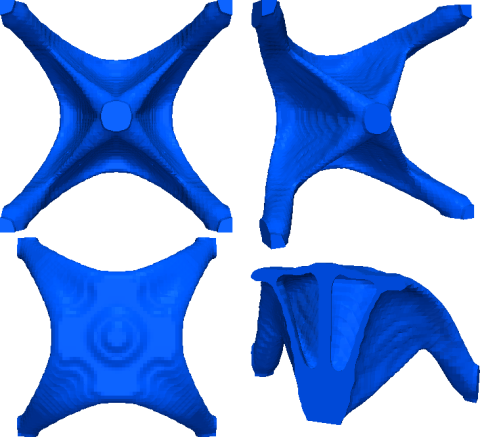}} \\
	\end{tabularx}
	\caption{Cube with center load SIMP topology designs. Clockwise: Bottom, Side, Top, and Clip views.}
	\label{fig:cube-center-load-views-SIMP}
\end{figure*}

The design domain is initialized with a uniform material distribution of $\rho_i = 0.1$. The optimized material distributions are shown in Figure \ref{fig:cube-center-load-views-SIMP} where material with a density lower than $\rho_i < 0.75$ is considered void. The strain energies are reported in Tab.~\ref{tab:cube-center-load-optimization-table}. For both meshes the volume constraint is active in the converged designs. As expected, the optimized geometry is smoother and the strain energy is lower for the refined mesh.

\begin{figure}
	\centering
	\includegraphics[width=\linewidth]{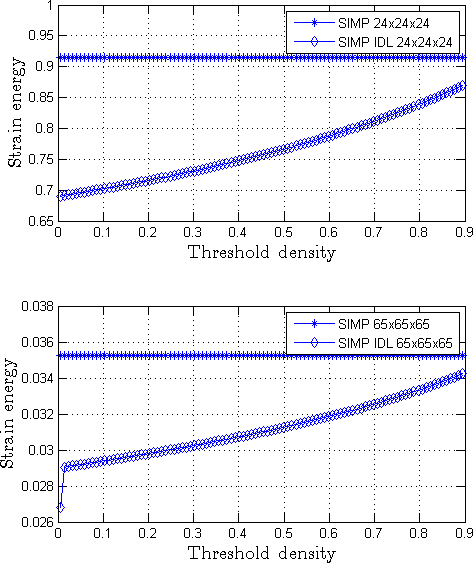}
	\caption{IDL post-processing for cube with center load problem.}
	\label{fig:cube-center-load-isovolumes-plot}
\end{figure}

The SIMP results for the coarse and fine mesh are post-processed with the IDL approach described above. The strain energies for varying threshold values, $\rho_T$, are plotted in Fig.~\ref{fig:cube-center-load-isovolumes-plot}. The volume constraint is met for $\rho_T = 0.78$ for the coarse mesh and $\rho_T = 0.44$ for the fine mesh. The difference between the SIMP-IDL and the SIMP results is of $0.04894 \%$ for the coarse mesh and $0.1230 \%$ for the fine mesh. The value of $\rho_T$ is higher for the coarse mesh because it cannot converge to a design with void inclusions. In both cases the strain energies of the post-processed results match well the SIMP predictions as the density distributions converged well to ``0-1'' solutions.

The volume fractions of intermediate densities are $0.285$ and $0.0189$ for the coarse and fine mesh, respectively. The post-processed designs have lower strain energies because the post-processing counteracts the effect of the density filter (\ref{eq:smoothing-equation-SIMP}).

\begin{figure}
	\centering
	\includegraphics[width=\linewidth]{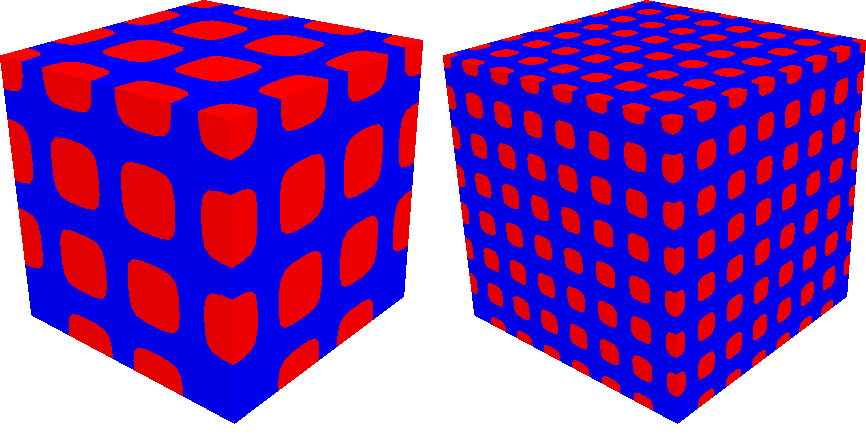}
	\caption{Initial level set configurations for the cube with a center load problem.}
	\label{fig:cube-center-load-initial-levelset}
\end{figure}

\begin{table}
	\centering
	\begin{tabular}{cc|c|c|c|c|c}
		\cline{2-6}
		& \multicolumn{1}{ |c| }{Mesh size}& \multicolumn{4}{ c| }{Strain energy}      \\ \cline{1-6}
		\multicolumn{1}{ |c| }{\multirow{2}{*}{SIMP} } &
		\multicolumn{1}{ |c| }{24x24x24} & \multicolumn{4}{ c| }{9.1456e-01} &     \\ \cline{2-6}
		\multicolumn{1}{ |c  }{}                        &
		\multicolumn{1}{ |c| }{65x65x65} & \multicolumn{4}{ c| }{3.5244e-02} &     \\ \cline{1-6}
		\multicolumn{1}{ |c  }{\multirow{2}{*}{XFEM} } &
		\multicolumn{1}{ |c| }{24x24x24} & \multicolumn{4}{ c| }{1.0082e+00} &     \\ \cline{2-6}
		\multicolumn{1}{ |c  }{}                        &
		\multicolumn{1}{ |c| }{65x65x65} & \multicolumn{4}{ c| }{3.5519e-02} &     \\ \cline{1-6}
	\end{tabular}
	\caption{Comparison of strain energy for the cube with a center load optimization problem for SIMP and LSM-XFEM approaches.}
\label{tab:cube-center-load-optimization-table}
\end{table}

\begin{figure*}
	\begin{tabularx}{\linewidth}{XX}
		\subfloat[24x24x24 mesh]{\includegraphics[width=\linewidth]{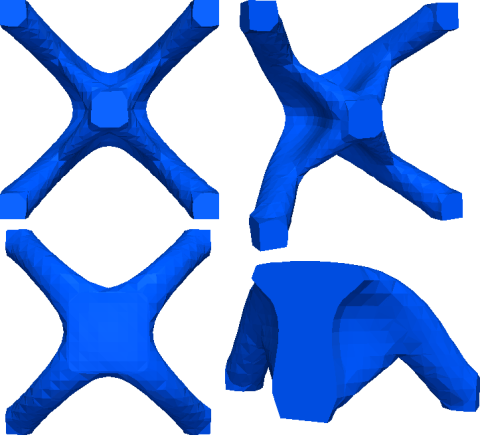}} &
		\subfloat[65x65x65 mesh]{\includegraphics[width=\linewidth]{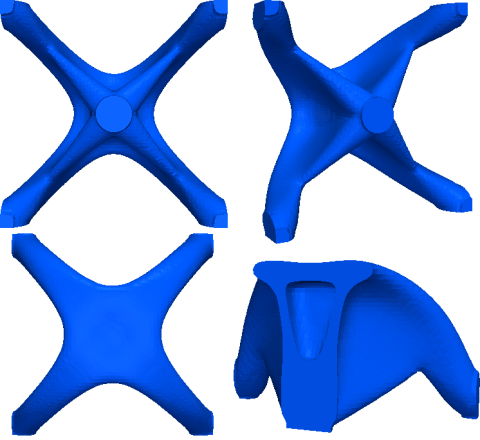}} \\
	\end{tabularx}
	\caption{Cube with center load XFEM topology designs. Clockwise: Bottom, Side, Top, and Clip views.}
	\label{fig:cube-center-load-views-XFEM}
\end{figure*}

\begin{figure}
	\centering
	\includegraphics[width=\linewidth]{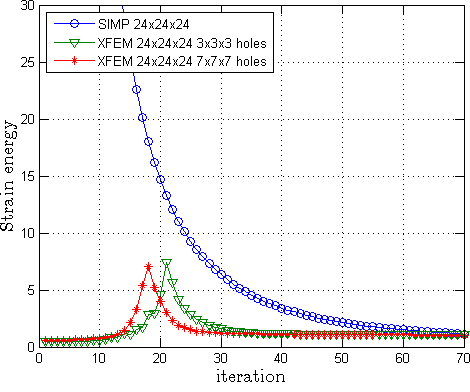}
	\caption{Convergence plot for the cube with a center load problem.}
	\label{fig:cube-center-load-XFEM-optimization-holes-coarse}
\end{figure}

The same optimization problem is solved with the proposed LSM-XFEM approach. The smoothing radius is set to $r_\phi = 3.2$ for the coarse mesh and $r_\phi = 1.182$ for the fine mesh. No perimeter constraint is imposed. We seed the initial design with two different configurations of void inclusions to study the influence of the initial layout on the optimization results. For both configurations we impose an equally spaced array of square-shaped holes with rounded corners, by modifying Eq. \ref{eqn:circle} into the following:

\begin{equation}
	\label{eqn:squircle}
	\phi_i = { (x_i - x_c) }^{10} + { (y_i - y_c) }^{10} + { (z_i - z_c) }^{10} - r^{10}.
\end{equation}

One configuration has $3 \times 3 \times 3$ equally spaced holes with radius $5.50$, the other $7 \times 7 \times 7$ holes with radius $2.0$, as shown in Fig.~\ref{fig:cube-center-load-initial-levelset}. In both cases, the volume constraint is not satisfied with the initial design. Note that no inclusions are placed at the four bottom corners where the boundary conditions are applied.

Both level set configurations converge to nearly indistinguishable designs and strain energy values, for both the coarse and fine meshes. The optimized designs are shown in Fig.~\ref{fig:cube-center-load-views-XFEM}. The strain energies of the optimized designs are given in Tab.~\ref{tab:cube-center-load-optimization-table}. The convergence history for the coarse meshes in SIMP and LSM-XFEM is shown in Fig.~\ref{fig:cube-center-load-XFEM-optimization-holes-coarse}.

For the example consider here, the SIMP and LSM-XFEM results match well, both in regards to the geometry and the strain energy values. The LSM-XFEM approach shows a faster convergence as the mesh is refined, for both the coarse and fine meshes. Comparing the optimized geometries, the SIMP results contain more structural features for both mesh resolutions. For example, considering the fine mesh, SIMP generates two small holes in the webs connecting the supports to the load point, while XFEM leads to only one larger hole, independent of the initial design configuration. However, these small differences have only a minor impact on the structural performance, i.e.~the strain energy, of the optimized designs.

Considering the conceptual structural layout, both, the SIMP and the LSM-XFEM approach, display only minor mesh dependencies for the problem studied here. Although the strain values show significant differences, the optimized geometries obtained with the coarse and fine meshes differ insignificantly for the SIMP and LSM-XFEM approach. The following example will demonstrate a less benign convergence and identify more pronounced differences between SIMP and LSM-XFEM.

%%%%%%%%%%%%%%%%%%%%%%%%%%%%%%%%%
% Building under torsion
%%%%%%%%%%%%%%%%%%%%%%%%%%%%%%%%%

\subsection{Cuboid under torsion}
\label{sec:building-under-torsion}

The second ``solid-void'' example is taken from \citet{NPS+:12} and reveals differences in the SIMP and LSM-XFEM approaches. We will show that, without imposing a mesh-independent minimum feature size constraint, the proposed LSM-XFEM approach may feature a significantly greater convergence than the SIMP method employed in this paper. However, we will also illustrate that our LSM-XFEM approach suffers from a lack of a robust and intuitive shape control technique.

\begin{figure}
	\centering
	\includegraphics[width=\linewidth]{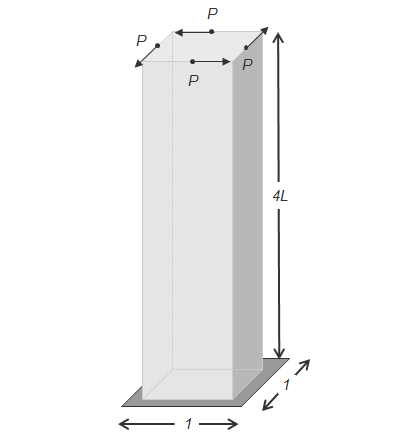}
	\caption{Cuboid under torsion model.}
	\label{fig:building-under-torsion-initial-setup}
\end{figure}

The design domain is a cuboid of size $4 \times 1 \times 1$, as shown in Fig.~\ref{fig:building-under-torsion-initial-setup}. A torque moment is generated via $4$ unit loads acting at the centers of the edges of the top face. The design domain is clamped at the bottom face. The Young's modulus is set to $1.0$ and the Poisson ratio to $0.3$. The volume of the stiff phase is constrained to $10$\% of the total volume. The problem is solved on the full mesh.

% Mesh convergence study %%%%%%%%

\subsubsection{Mesh convergence study}
\label{sec:mesh-convergence-study}

The optimization problem is solved with the SIMP approach for four different mesh sizes: $40 \times 10 \times 10$, $60 \times 15 \times 15$,  $80 \times 20 \times 20$, and $120 \times 30 \times 30$. The Young's modulus of the void phase is set to $E_B = 10^{-9}$. The design domain is initialized with a uniform material distribution of $\rho_i = 0.1$. The penalization factor is $p=3$. First we consider a projection parameter of $\beta=0$ and scale the smoothing radius with the element size: $r_\rho = 1.6 \times$ the element edge length.

\begin{figure}
	\begin{tabularx}{\linewidth}{XX}
		\subfloat[ 40x10x10]{\includegraphics[width=\linewidth]{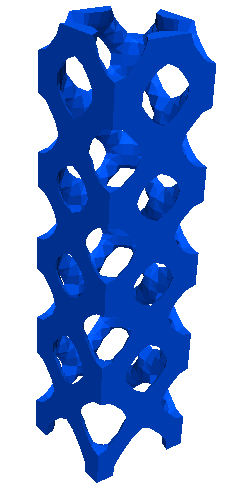}} &
		\subfloat[ 60x15x15]{\includegraphics[width=\linewidth]{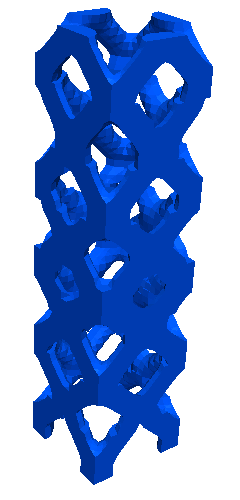}} \\
		\subfloat[ 80x20x20]{\includegraphics[width=\linewidth]{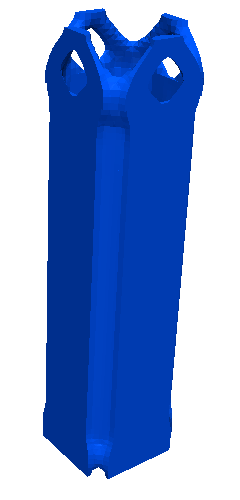}} &
		\subfloat[120x30x30]{\includegraphics[width=\linewidth]{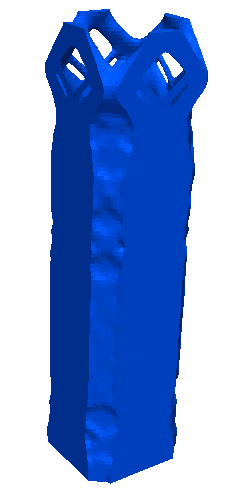}} \\
	\end{tabularx}
	\caption{Cuboid under torsion results with SIMP. $p=3$, $r_\rho = 1.6$ of the element edge length, $\beta=0$.}
	\label{fig:building-under-torsion-SIMP-optimization-meshDependent}
\end{figure}

The optimized material distributions are shown in Fig.~\ref{fig:building-under-torsion-SIMP-optimization-meshDependent} where material with a density lower than $\rho_i < 0.35$ is considered void. The strain energies are reported in Tab.~\ref{tab:building-under-torsion-strain-energy-SIMP-mesh-dependent} and display the expected decrease in strain energy as the mesh is refined. For all meshes the volume constraint is active in the converged designs.

\begin{table}
	\centering
	\begin{tabular}{cc|c|c|c|c|c}
		\cline{2-6}
		& \multicolumn{1}{ |c| }{Mesh size}& \multicolumn{4}{ c| }{Strain energy}       \\ \cline{1-6}
		\multicolumn{1}{ |c| }{\multirow{4}{*}{SIMP} } &
		\multicolumn{1}{ |c| }{ 40x10x10} & \multicolumn{4}{ c| }{7.5195e+03} &     \\ \cline{2-6}
		\multicolumn{1}{ |c  }{}                        &
		\multicolumn{1}{ |c| }{ 60x15x15} & \multicolumn{4}{ c| }{4.2076e+03} &     \\ \cline{2-6}
		\multicolumn{1}{ |c  }{}                        &
		\multicolumn{1}{ |c| }{ 80x20x20} & \multicolumn{4}{ c| }{4.0298e+03} &     \\ \cline{2-6}
		\multicolumn{1}{ |c  }{}                        &
		\multicolumn{1}{ |c| }{120x30x30} & \multicolumn{4}{ c| }{2.6555e+03} &     \\ \cline{1-6}
	\end{tabular}
	\caption{Strain energy results for the cuboid under torsion SIMP problem with $p=3$, $r_\rho = 1.6$ of the element edge length, and $\beta=0$.}
	\label{tab:building-under-torsion-strain-energy-SIMP-mesh-dependent}
\end{table}

As the mesh is refined, the evolution of the SIMP results shows an interesting discontinuity which is typically not observed for two dimensional problems. The optimized material layout switches abruptly from a grid-type structure, which conceptually agrees with the results of \citet{NPS+:12}, to a hollow square prism design. In contrast to two dimensional structures, where refining the mesh with a mesh-dependent filter radius leads to an ever increasing number of holes, in this example the opposite is the case. As the element size drops below a threshold, it is more advantageous to form a continuous thin outer wall rather than a grid-type structure. This behavior is a direct consequence of the combination of SIMP penalization and density smoothing. We will revisit this issue again later.

\begin{figure}
	\begin{tabularx}{\linewidth}{XX}
		\subfloat[ 40x10x10]{\includegraphics[width=\linewidth]{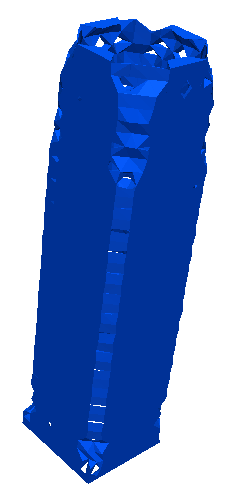}} &
		\subfloat[120x30x30]{\includegraphics[width=\linewidth]{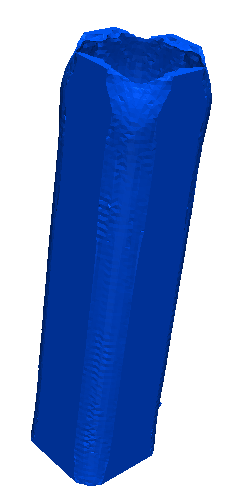}} \\
	\end{tabularx}
	\caption{Cuboid under torsion results with the LSM-XFEM approach. $r_\phi = 1.6$ of the element edge length, no perimeter constraint.}
	\label{fig:building-under-torsion-XFEM-noPerimeter-meshDependent}
\end{figure}

The LSM-XFEM results for a smoothing radius of $r_\phi = 1.6 \times$ the element edge length are shown in Fig.~\ref{fig:building-under-torsion-XFEM-noPerimeter-meshDependent}. No perimeter constraint is applied to this problem.  Here only the results for the coarsest and the finest meshes of the SIMP study above are shown. Note that in contrast to the SIMP results, the LSM-XFEM approach leads to conceptually equivalent design on both meshes. Refining the mesh only improves some local details. This feature is due to the ability of the LSM to represent thin structural features on coarse meshes. The thickness of the walls is measured by manually comparing the coordinates of two points at half the height of the design in the visualization tool. The thickness value is $0.0288$ for the coarse mesh and $0.0276$ for the fine mesh.

The strain energies of the LSM-XFEM results are given in Tab.~\ref{tab:building-under-torsion-strain-energy-XFEM-mesh-dependent}. The strain energy for the fine mesh is slightly larger than the one of the coarse mesh. This effect is due to the tendency of coarse finite element discretization over predicting the stiffness.

\begin{table}
	\centering
	\begin{tabular}{cc|c|c|c|c|c}
		\cline{2-6}
		& \multicolumn{1}{ |c| }{Mesh size}& \multicolumn{4}{ c| }{Strain energy}       \\ \cline{1-6}
		\multicolumn{1}{ |c  }{\multirow{2}{*}{XFEM} } &
		\multicolumn{1}{ |c| }{ 40x10x10} & \multicolumn{4}{ c| }{8.7551e+02} &     \\ \cline{2-6}
		\multicolumn{1}{ |c  }{}                        &
		\multicolumn{1}{ |c| }{120x30x30} & \multicolumn{4}{ c| }{9.8262e+02} &     \\ \cline{1-6}
	\end{tabular}
	\caption{Strain energy results for cuboid under torsion XFEM problem with $r_\phi = 1.6$ of the element
	edge length, no perimeter constraint.}
\label{tab:building-under-torsion-strain-energy-XFEM-mesh-dependent}
\end{table}

The differences between the SIMP and LSM-XFEM results are significant. Although the discrepancy in strain energy decreases as the mesh is refined, the difference is large even for the two finer meshes where the SIMP and LSM-XFEM results are similar. As the following investigation will show, the poorer performance of the SIMP results is primarily caused by the density filter, which prevents the material distribution to converge to a ``0-1'' result.

\begin{figure*}
	\begin{tabularx}{\linewidth}{XXX}
		\subfloat[$\beta = 0$]{\includegraphics[width=\linewidth]{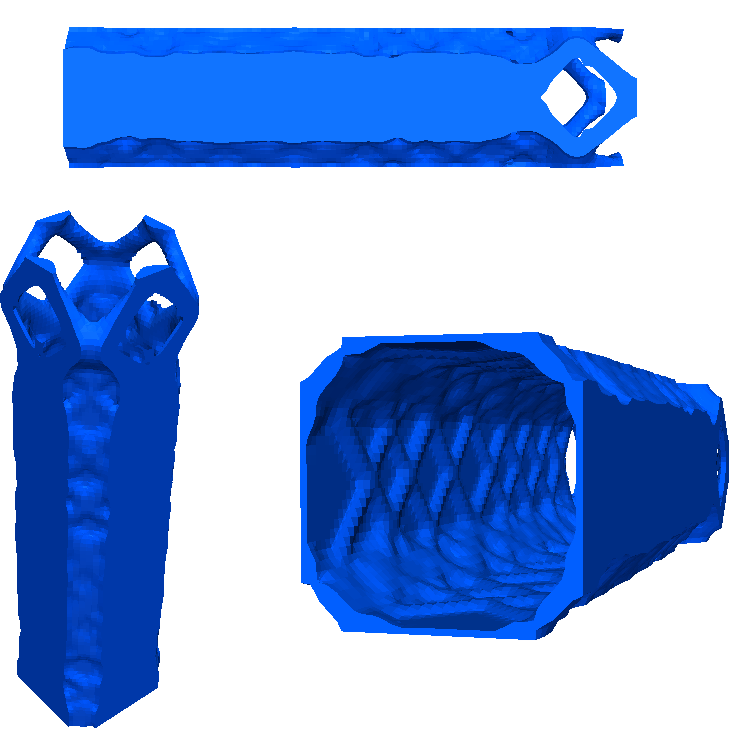}} &
		\subfloat[$\beta = 4$]{\includegraphics[width=\linewidth]{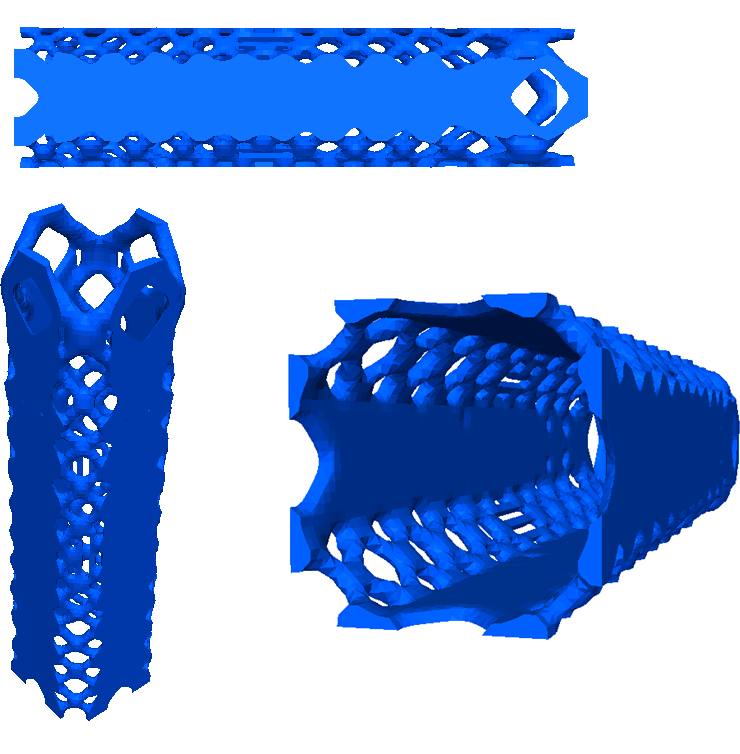}} &
		\subfloat[$\beta = 8$]{\includegraphics[width=\linewidth]{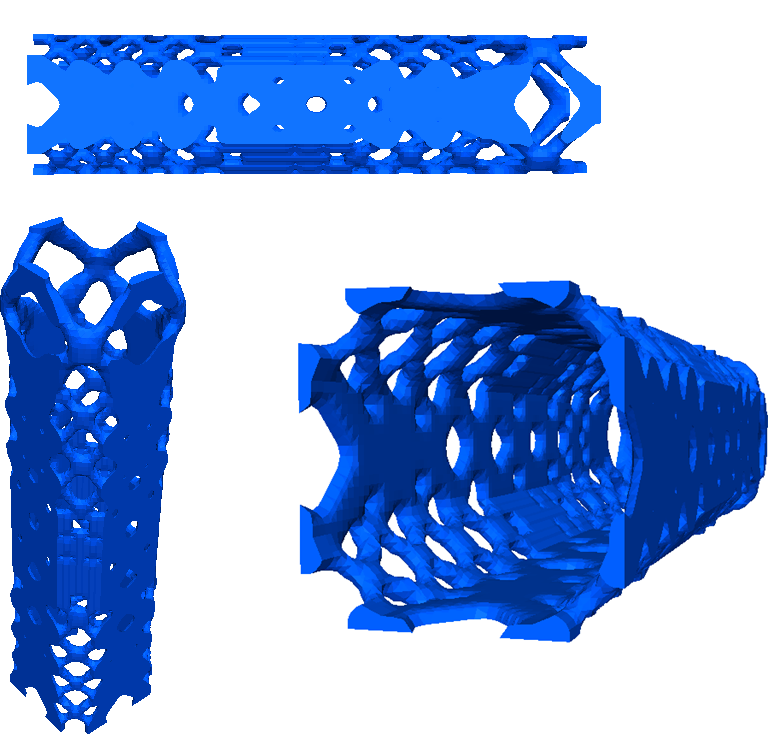}} \\
	\end{tabularx}
	\caption{Density isovolumes of 0.35 for the three step SIMP results.}
	\label{fig:building-under-torsion-projection}
\end{figure*}

First we study the influence of the projection scheme (\ref{eq:heaviside-SIMP}) on the SIMP results for the most refined mesh. The optimized material distributions for $\beta=4.0$ and $\beta=8.0$ are shown in Fig.~\ref{fig:building-under-torsion-projection}, where material with a density lower than $\rho_i < 0.35$ is considered void. For convenience the result for $\beta=0.0$ is shown again. Table \ref{tab:building-under-torsion-strain-energy-projection} reports on the strain energies and the volume fractions of intermediate densities, $\bar{\rho}$, as the projection parameter, $\beta$, is increased. The higher $\beta$, the lower $\bar{\rho}$ and the lower the strain energy, approaching the one of the LSM-XFEM result. Note that as $\beta$ increases, the more holes emerge. The thickness of the walls for $\beta=8$ is $0.0434$, which is smaller than the value for $\beta = 0$, $0.0447$, and closer to the XFEM value.

\begin{table}
	\centering
	\begin{tabular}{cc|c|c|c|c|c|c|c|c}
		\cline{2-9}
		& \multicolumn{1}{ |c| }{$\beta$ projection} & \multicolumn{3}{ c| }{strain energy} & \multicolumn{4}{ c| }{$\bar{\rho}$ utilization} \\ \cline{1-9}
		\multicolumn{1}{ |c| }{\multirow{3}{*}{SIMP} } &
		\multicolumn{1}{ |c| }{0} & \multicolumn{3}{ c| }{2.6555e+03} & \multicolumn{4}{ c| }{4.0688e-02} & \\ \cline{2-9}
		\multicolumn{1}{ |c  }{}                        &
		\multicolumn{1}{ |c| }{4} & \multicolumn{3}{ c| }{2.0264e+03} & \multicolumn{4}{ c| }{2.5131e-02} & \\ \cline{2-9}
		\multicolumn{1}{ |c  }{}                        &
		\multicolumn{1}{ |c| }{8} & \multicolumn{3}{ c| }{1.9039e+03} & \multicolumn{4}{ c| }{1.9509e-02} & \\ \cline{1-9}
	\end{tabular}
	\caption{Strain energy results for the cuboid under torsion SIMP problem with the projection scheme of Eq. \ref{eq:heaviside-SIMP}.}
\label{tab:building-under-torsion-strain-energy-projection}
\end{table}

\begin{figure}
	\centering
	\includegraphics[width=\linewidth]{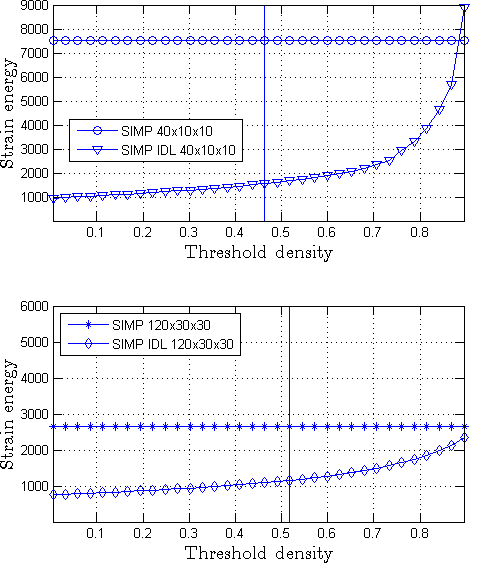}
	\caption{SIMP-IDL post-processing for the building under torsion problem. The vertical line marks $\rho_T = 0.46$, at
	which the volume constraint is satisfied for both meshes.}
	\label{fig:building-under-torsion-isovolumes-plot}
\end{figure}

Instead of enforcing a better convergence toward a ``0-1'' solution by increasing the projection parameter $\beta$, we post-process the SIMP results for $\beta=0$ by the IDL post-processing method. Figure \ref{fig:building-under-torsion-isovolumes-plot} shows the strain energy of the post-processed design for the coarsest and the finest mesh. The volume constraint is satisfied for a threshold value of $\rho_T=0.4634$ for the coarse mesh, and $\rho_T=0.5174$ for the fine mesh. The difference between the SIMP-IDL and the SIMP results is of $0.7066 \%$ for the coarse mesh and $0.6282 \%$ for the fine mesh at the threshold $\rho_T$. The associated strain energies are given in Tab.~\ref{tab:building-under-torsion-strain-energy-SIMP-IDL}.

\begin{table}
	\centering
	\begin{tabular}{cc|c|c|c|c|c|c|c|c}
		\cline{2-9}
		& \multicolumn{1}{ |c| }{mesh} & \multicolumn{3}{ c| }{$\rho_T$} & \multicolumn{4}{ c| }{strain energy} \\
		\cline{1-9} \multicolumn{1}{ |c| }{\multirow{2}{*}{SIMP} } &
		\multicolumn{1}{ |c| }{ 40x10x10} & \multicolumn{3}{ c| }{0.4634} & \multicolumn{4}{ c| }{1.5601e+03} & \\
		\cline{2-9} \multicolumn{1}{ |c  }{}                        &
		\multicolumn{1}{ |c| }{120x30x30} & \multicolumn{3}{ c| }{0.5174} & \multicolumn{4}{ c| }{1.1530e+03} & \\
		\cline{1-9}
	\end{tabular}
	\caption{Strain energy at $\rho_T = 0.4634$ and $\rho_T = 0.5174$ for the cuboid under torsion SIMP problem.}
\label{tab:building-under-torsion-strain-energy-SIMP-IDL}
\end{table}

The strain energy of the post-processed results of the fine mesh is rather similar to the result obtained for SIMP with $\beta=8.0$ in Tab.~\ref{tab:building-under-torsion-strain-energy-projection} and the LSM-XFEM results in Tab.~\ref{tab:building-under-torsion-strain-energy-XFEM-mesh-dependent}. For the coarse mesh, the strain energy is well below the raw SIMP results from Tab.~\ref{tab:building-under-torsion-strain-energy-SIMP-mesh-dependent} but still above the results for the LSM-XFEM approach in Tab.~\ref{tab:building-under-torsion-strain-energy-XFEM-mesh-dependent}. As we will see below, this is because of the larger smoothing radius which prevents the formation of smaller features and thinner walls.

% Feature size control %%%%%%%%%%

\subsubsection{Feature size control}
\label{sec:feature-size-control}

The mesh refinement study above suggests that the results of the LSM-XFEM approach are less sensitive to mesh refinement than the SIMP method without mesh-independent filtering. Geometric features, such as the thin walls, can be represented on coarse and fine meshes, independent of their size. This observation is in agreement with studies for two-dimensional problems, see for example \citet{KM:12}, but the phenomena is more pronounced and of greater importance for three dimensional problems. While this feature of the LSM-XFEM approach may be considered an advantage, the ability to control the minimum feature size is of importance for many applications and to account for manufacturing constraints and costs. The following study will show that the proposed LSM-XFEM approach currently lacks the ability to efficiently and intuitively control the local feature size.

\begin{figure}
	\begin{tabularx}{\linewidth}{XX}
		\subfloat[SIMP]{\includegraphics[width=\linewidth]{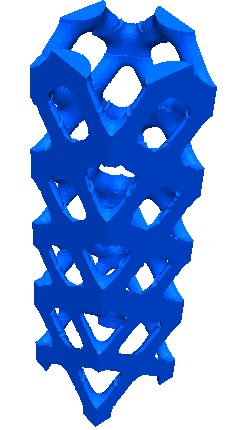}} &
		\subfloat[XFEM]{\includegraphics[width=\linewidth]{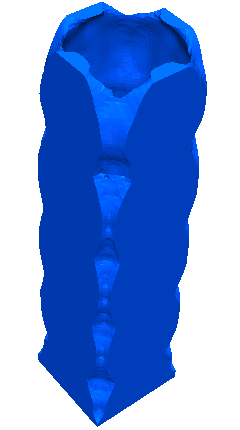}} \\
	\end{tabularx}
	\caption{Topology design with a larger smoothing radius.}
	\label{fig:building-under-torsion-radius}
\end{figure}

We first show that applying the same absolute filter radius in the SIMP formulation efficiently controls the feature size. Figure \ref{fig:building-under-torsion-radius}(a) shows the SIMP results on the $120 \times 30 \times 30$ mesh for a projection parameter $\beta=0$, a penalization factor of $p=3$, and a smoothing radius of $r_\rho = 0.16$ which is the same radius applied earlier for the coarsest mesh in Fig.~\ref{fig:building-under-torsion-SIMP-optimization-meshDependent}(a). Comparing the SIMP results in Fig.~\ref{fig:building-under-torsion-radius}(a) and Fig.~\ref{fig:building-under-torsion-SIMP-optimization-meshDependent}(a) confirms the finding of numerous studies \citep{BS:03} that the SIMP approach leads to the same conceptual layout independent of the mesh refinement level if a mesh-independent filter is used. The strain energy of the design in Fig.~\ref{fig:building-under-torsion-radius}(a) is given in Tab.~\ref{tab:building-under-torsion-radius}.

\begin{table}
	\centering
	\begin{tabular}{cc|c|c|c|c|c}
		\cline{2-6}
		& \multicolumn{1}{ |c| }{Mesh size}& \multicolumn{4}{ c| }{Strain energy}      \\ \cline{1-6}
		\multicolumn{1}{ |c| }{\multirow{1}{*}{SIMP} } &
		\multicolumn{1}{ |c| }{120x30x30} & \multicolumn{4}{ c| }{6.5772e+03} &     \\ \cline{1-6}
		\multicolumn{1}{ |c  }{\multirow{1}{*}{XFEM} } &
		\multicolumn{1}{ |c| }{120x30x30} & \multicolumn{4}{ c| }{8.2077e+02} &     \\ \cline{1-6}
	\end{tabular}
	\caption{Strain energy for the cuboid under torsion optimization problem with a mesh independent filter.}
\label{tab:building-under-torsion-radius}
\end{table}

A similar effect is not observed in the LSM-XFEM approach when we apply the same filter radius, $r_\phi$, used for the coarse mesh to the fine mesh. Figure \ref{fig:building-under-torsion-radius}(b) shows the outcome of this procedure. The overall design is unchanged, and increasing the smoothing radius results in a less smooth design. The strain energy of this design is reported in Tab.~\ref{tab:building-under-torsion-radius}.

\begin{figure}
	\begin{tabularx}{\linewidth}{XX}
		\subfloat[]{\includegraphics[width=\linewidth]{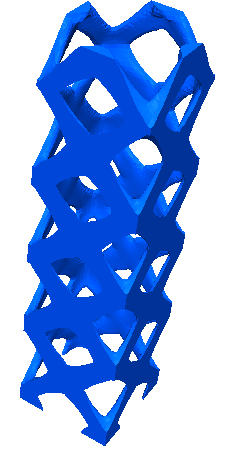}} &
		\subfloat[]{\includegraphics[width=\linewidth]{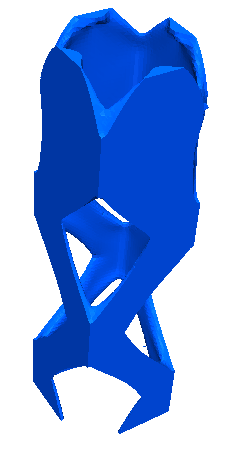}} \\
	\end{tabularx}
	\caption{Cuboid under torsion problem with a perimeter constraint: (a) XFEM restarted from SIMP and (b) XFEM restarted from XFEM.}
	\label{fig:building-under-torsion-XFEM-perimeter}
\end{figure}

To control the overall structural complexity in the LSM, the formulation of the optimization problem (\ref{eq:topo_problem}) is often augmented by a perimeter constraint \citep{DML+:13}. While this approach does not directly control the minimum feature size, reducing the maximum feasible perimeter often removes small features which do not alter much the structural performance. To study the influence of a perimeter constraint on the torsion problem, we perform the following two numerical experiments on the $120 \times 30 \times 30$ mesh using the LSM-XFEM approach. We measure the perimeter of the SIMP result shown in Fig.~\ref{fig:building-under-torsion-radius}(a) and impose this value as an upper bound on the perimeter. One problem uses the SIMP result in Fig.~\ref{fig:building-under-torsion-radius}(a) as the initial design, and the other one uses the LSM-XFEM result in Fig.~\ref{fig:building-under-torsion-XFEM-noPerimeter-meshDependent}(b). The results are shown in Fig.~\ref{fig:building-under-torsion-XFEM-perimeter} and the strain energies are given in Tab.~\ref{tab:building-under-torsion-XFEM-perimeter}.

Depending on the initial designs, the LSM-XFEM problems converge to different designs. While the design in Fig.~\ref{fig:building-under-torsion-XFEM-perimeter}(b) displays a truss-like design in the bottom half of the design domain, the perimeter constraint does not prevent the formation of thin walls in the upper half. The thickness of the walls in the upper half of the design is $0.0188$. Thus, the perimeter constraint does not control the local feature size. The design in Fig.~\ref{fig:building-under-torsion-XFEM-perimeter}(a) resembles closely the SIMP result from which it was restarted. However, considering the strain energy in Tab.~\ref{tab:building-under-torsion-XFEM-perimeter}, this design has a larger strain energy than the one in Fig.~\ref{fig:building-under-torsion-XFEM-perimeter}(b).

The study above has shown that neither smoothing the level set field nor imposing a perimeter constraint allows controlling the minimum feature size. Further, the effect of a perimeter constraint is non-intuitive as the result in Fig.~\ref{fig:building-under-torsion-XFEM-perimeter}(b) shows. The design has more structural features than the design without perimeter constraint in Fig.~\ref{fig:building-under-torsion-XFEM-noPerimeter-meshDependent}(b).

\begin{table}
	\centering
	\begin{tabular}{cc|c|c|c|c|c}
		\cline{2-6}
		& \multicolumn{1}{ |c| }{Initial design}& \multicolumn{4}{ c| }{Strain energy}      \\ \cline{1-6}
		\multicolumn{1}{ |c| }{\multirow{2}{*}{XFEM} } &
		\multicolumn{1}{ |c| }{SIMP} & \multicolumn{4}{ c| }{1.3321e+03} &     \\ \cline{2-6}
		\multicolumn{1}{ |c  }{} &
		\multicolumn{1}{ |c| }{XFEM} & \multicolumn{4}{ c| }{1.3185e+03} &     \\ \cline{1-6}
	\end{tabular}
	\caption{Strain energy for the cuboid under torsion optimization problem with XFEM and a perimeter constraint.}
\label{tab:building-under-torsion-XFEM-perimeter}
\end{table}

%%%%%%%%%%%%%%%%%%%%%%%%%%%%%%%%%
% Two phase optimization
%%%%%%%%%%%%%%%%%%%%%%%%%%%%%%%%%

\subsection{Two-phase Cantilevered Beam Design}
\label{sec:two-phase-optimization}

\begin{figure}
	\centering
	\includegraphics[width=\linewidth]{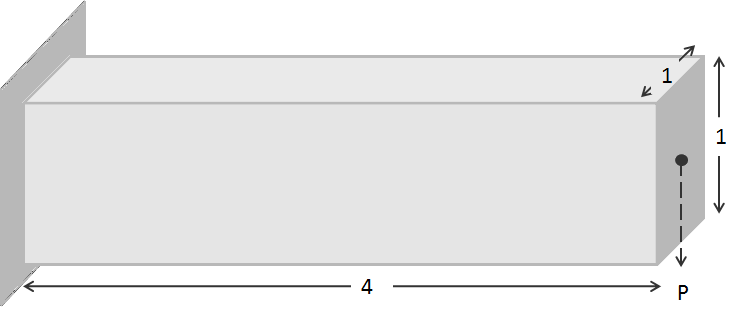}
	\caption{Initial setup for the long cantilever beam two-phase problem.}
	\label{fig:long-cantilever-beam-initial-setup}
\end{figure}

The examples in the two previous subsections were concerned with ``solid-void'' problems. Here we study a ``solid-solid'' problem to demonstrate the applicability of the proposed LSM-XFEM approach to this class of problems. Note that the simplified XFEM formulation discussed in Section \ref{sec:intro} is not applicable to such problems. The generalized enrichment strategy of Section \ref{sec:discretization} is required and the interface conditions of Eq. \ref{eq:stabilized_Residual} need to be satisfied.

\begin{figure*}[t]
	\centering
	\includegraphics[width=\linewidth]{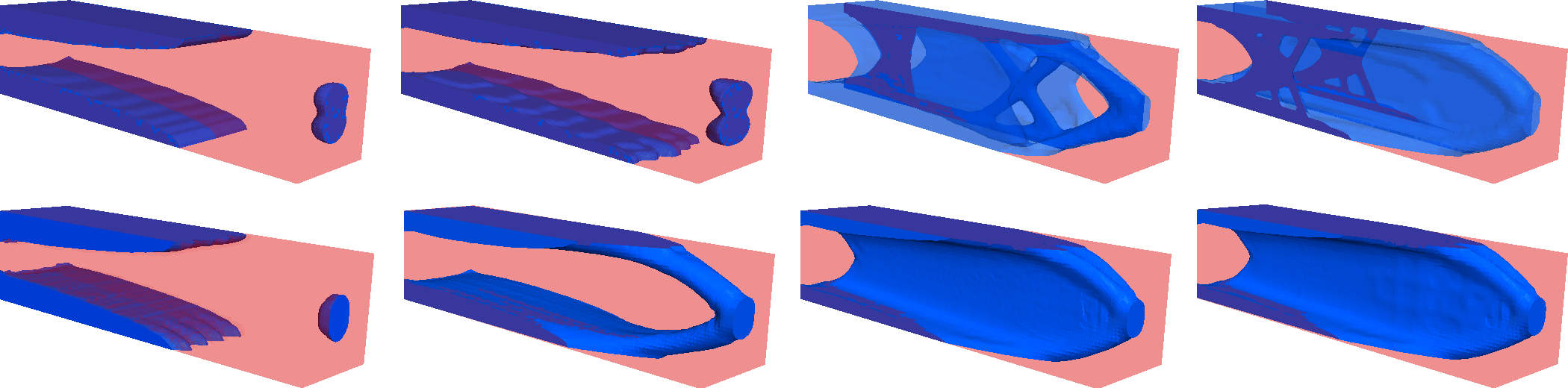}
	\caption{Two-phase topology optimization using SIMP and XFEM. First row shows the SIMP results, second row, XFEM.
	Density ratios are from left to right $E_{B} = 0.5 E_{A}$, $E_{B} = 0.1 E_{A}$, $E_{B} = 0.01 E_{A}$, $E_{B}$ is void.}
	\label{fig:long-cantilever-beam-SIMP-XFEM}
\end{figure*}

We study the optimal two-phase layout of a $4 \times 1 \times 1$ cantilevered beam subject to a tip load; see Fig.~\ref{fig:long-cantilever-beam-initial-setup}. The stiff  phase ``A'' has Young's modulus of $E_A = 1.0$; three values of Young's moduli for the soft phase are considered:  $E_B = [0.5, 0.1, 0.01]$. Both phases have a Poisson ratio of $0.3$. The maximum volume of the stiff phase is limited to $30 \%$ of the total volume. The design domain is discretized by $120 \times 30 \times 30$ elements. Because of the symmetry condition, only one half of the cuboid is numerically analyzed. We compare the SIMP and LSM-XFEM results.

The optimization problem is solved by a SIMP approach with a penalization factor of $p=3$, a smoothing radius of $r_{\rho} = 0.05333$ ($1.6 \times$ the element edge) and the projection parameter of $\beta=0$. The design domain is initialized with a uniform material distribution of $\rho_i = 0.3$.  The optimized material distributions are shown in Fig.~\ref{fig:long-cantilever-beam-SIMP-XFEM} where material with a density lower than $\rho_i < 0.25$ is transparent. The strain energies are reported in Tab.~\ref{tab:long-cantilever-beam-two-phase-optimization}.

The LSM-XFEM results for a smoothing radius of $r_\phi = 0.05333$ ($1.6 \times$ the element edge) are shown in Fig.~\ref{fig:long-cantilever-beam-SIMP-XFEM} and the strain energies are given in Tab.
\ref{tab:long-cantilever-beam-two-phase-optimization}. Considering the full design domain, the level set field is initialized with a $16 \times 4 \times 4$ array of equally spaced holes with radius of $0.1050$. The initial design is shown in Fig.~\ref{fig:long-cantilever-beam-swiss-cheese} and satisfies the volume constraint for the stiff phase. Note that the interface condition is enforced via the stabilized Lagrange multiplier method (\ref{eq:stabilized_Residual}) with an element wise constant Lagrange multiplier, $\lambda$, and a consistency factor of $\gamma = 10 \ (E_A+E_B)$.

Comparing the SIMP and LSM-XFEM results, the same trends can be observed for this ``solid-solid'' problem as for the ``solid-void'' ones studied earlier. The LSM-XFEM approach leads to three dimensional structures with thinner walls and higher stiffness. In contrast, the SIMP method generates truss-type structures, in particular if the discretization is too coarse and the optimum wall thickness is less than the size of an element.

\begin{figure*}
	\centering
	\includegraphics[width=0.6\linewidth]{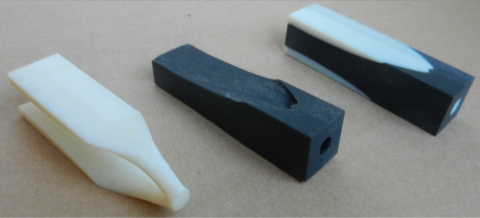}
	\caption{LSM-XFEM optimized two-phase design for $E_{B} = 0.1 E_{A}$ realized by 3D printing.}
	\label{fig:manufacture-long-cantilever-beam}
\end{figure*}

For illustration purposes only, we show a realization of the LSM-XFEM optimized design for $E_{B} = 0.1 E_{A}$ in Fig.~\ref{fig:manufacture-long-cantilever-beam}. The structure was fabricated with a polyjet 3D printing process on a Connex Objet 260 printer. White material represents phase ``A'', black represents phase ``B''. The left and center pieces show the individual phases printed separately, the printed two-phase design is shown on the right. 

\begin{table}
	\centering
	\begin{tabular}{cc|c|c|c|c|c}
		\cline{2-6}
		& \multicolumn{1}{ |c| }{Density ratio}& \multicolumn{4}{ c| }{Strain energy}       \\ \cline{1-6}
		\multicolumn{1}{ |c| }{\multirow{4}{*}{SIMP} } &
		\multicolumn{1}{ |c| }{$E_{B} = 0.50 E_{A}$} & \multicolumn{4}{ c| }{4.4081e-05} &     \\ \cline{2-6}
		\multicolumn{1}{ |c  }{}                        &
		\multicolumn{1}{ |c| }{$E_{B} = 0.10 E_{A}$} & \multicolumn{4}{ c| }{6.2862e-05} &     \\ \cline{2-6}
		\multicolumn{1}{ |c  }{}                        &
		\multicolumn{1}{ |c| }{$E_{B} = 0.01 E_{A}$} & \multicolumn{4}{ c| }{7.8627e-05} &     \\ \cline{2-6}
		\multicolumn{1}{ |c  }{}                        &
		\multicolumn{1}{ |c| }{$E_{B}$ is void     } & \multicolumn{4}{ c| }{7.6721e-05} &     \\ \cline{1-6}
		\multicolumn{1}{ |c  }{\multirow{4}{*}{XFEM} } &
		\multicolumn{1}{ |c| }{$E_{B} = 0.50 E_{A}$} & \multicolumn{4}{ c| }{4.3221e-05} &     \\ \cline{2-6}
		\multicolumn{1}{ |c  }{}                        &
		\multicolumn{1}{ |c| }{$E_{B} = 0.10 E_{A}$} & \multicolumn{4}{ c| }{5.9192e-05} &     \\ \cline{2-6}
		\multicolumn{1}{ |c  }{}                        &
		\multicolumn{1}{ |c| }{$E_{B} = 0.01 E_{A}$} & \multicolumn{4}{ c| }{6.4448e-05} &     \\ \cline{2-6}
		\multicolumn{1}{ |c  }{}                        &
		\multicolumn{1}{ |c| }{$E_{B}$ is void     } & \multicolumn{4}{ c| }{6.6283e-05} &     \\ \cline{1-6}
	\end{tabular}
	\caption{Two-phase topology optimization using SIMP and XFEM.}
\label{tab:long-cantilever-beam-two-phase-optimization}
\end{table}

%%%%%%%%%%%%%%%%%%%%%%%%%%%%%%%%%%%%%%%%%%%%%%%%%%%%%%%%%%%%%%%%%%%
% Conclusions
%%%%%%%%%%%%%%%%%%%%%%%%%%%%%%%%%%%%%%%%%%%%%%%%%%%%%%%%%%%%%%%%%%%

\section{Conclusions}
\label{sec:conclusions}

In this paper we presented an optimization approach combining a level set method (LSM) for describing the geometry and an extended finite element method (XFEM) for predicting the structural response. Building upon generalized enrichment and preconditioning schemes, previously developed for two-dimensional problems, the proposed optimization scheme was applied to two-phase ``solid-void'' and ``solid-solid'' problems in three dimensions. In all examples, the strain energy was minimized subject to a volume constraint on the stiff phase. The results of the LSM-XFEM approach with and without perimeter constraints were compared with the ones of a SIMP method which employs density filtering and projection.

The numerical studies suggest that the LSM-XFEM method features an improved convergence as the mesh is refined and is able to represent thin-walled structures on coarse meshes. The SIMP approach may require a strong projection to achieve clear ``0-1'' results with comparable strain energies. While density filtering is an efficient and intuitive method to control the local feature size, neither level set smoothing nor imposing a perimeter constraint achieve a similar effect on LSM-XFEM results.

The current lack of a feature size control and the significant improved complexity of the LSM-XFEM formulation limit the attractiveness of this scheme. However, for problems where a high mesh resolution is not tolerable and/or interface conditions need to be enforced with high accuracy, the LSM-XFEM approach might be an interesting alternative to SIMP-type methods. The advantages of the LSM-XFEM problem have been shown by \citet{KM:12} for fluid problems at high Reynolds numbers in two dimensions. The authors plan to study three dimensional flow problems with the LSM-XFEM approach in the future.

%%%%%%%%%%%%%%%%%%%%%%%%%%%%%%%%%%%%%%%%%%%%%%%%%%%%%%%%%%%%%%%%%%%
% Bibliography
%%%%%%%%%%%%%%%%%%%%%%%%%%%%%%%%%%%%%%%%%%%%%%%%%%%%%%%%%%%%%%%%%%%

%\bibliographystyle{spbasic}
%\bibliography{JabRefDatabase}

\end{document}